\documentstyle{amsppt}
\magnification=1200
\hsize=150truemm
\vsize=224.4truemm
\hoffset=4.8truemm
\voffset=12truemm

\NoRunningHeads

\define\C{{\bold C}}
 
\define\R{{\bold R}}
\define\Z{{\bold Z}}
 
\let\thm\proclaim
\let\fthm\endproclaim

\newcount\tagno
\newcount\secno
\newcount\subsecno
\newcount\stno
\global\subsecno=1
\global\tagno=0
\define\ntag{\global\advance\tagno by 1\tag{\the\tagno}}

\define\sta{\ 
{\the\secno}.\the\stno
\global\advance\stno by 1}

\define\stas{\the\stno
\global\advance\stno by 1}

\define\sect{\global\advance\secno by 1
\global\subsecno=1\global\stno=1\
{\the\secno}. }

\def\nom#1{\edef#1{{\the\secno}.\the\stno}}
\def\inom#1{\edef#1{\the\stno}}
\def\eqnom#1{\edef#1{(\the\tagno)}}

\newcount\refno
\global\refno=0

\def\nextref#1{
      \global\advance\refno by 1
      \xdef#1{\the\refno}}

\def\bref {\ref\global\advance\refno by 1\key{\the\refno}}


\nextref\AGM
\nextref\DON
\nextref\DS
\nextref\GA
\nextref\HOR

\topmatter
\title 
Rational convexity of non generic immersed lagrangian submanifolds
 \endtitle
\author Julien Duval
 \footnote{Laboratoire de Math\'ematique, Universit\'e Paris-Sud, 91405 Orsay Cedex \newline julien.duval\@math.u-psud.fr}
 and Damien Gayet
\footnote{Institut Camille Jordan, Universit\'e Claude Bernard, 69622 Villeurbanne Cedex \newline gayet\@math.univ-lyon1.fr}
\footnote""{Keywords: rational convexity, lagrangian submanifold. AMS class.: 32E20, 53D12 \newline}
\endauthor
\abstract \  We prove that an immersed lagrangian submanifold
 in $\C^n$ with quadratic self-tangencies is rationally convex. This generalizes former results for the embedded and the immersed transversal cases.
\endabstract 
  \endtopmatter 
\document
 \subhead 0. Introduction \endsubhead
\stno=1

\null

\noindent

Let $K$ a compact set in $\C^n$. It is {\it rationally convex} if its complement
can be filled out by global complex hypersurfaces. According to the Oka-Weil
theorem, holomorphic functions nearby such a rationally convex compact set $K$ can be approximated
on $K$ by rational functions.
  
A classical obstruction to the rational convexity of $K$ is the presence in $\C^n$
of a Riemann surface $C$ with boundary in $K$, such that $\partial C$  bounds in $K$. This obstruction disappears
when $K$ is lagrangian for a global K\"ahler form, i.e. when $K$ is a submanifold of real dimension $n$
on which this form vanishes. Indeed N. Sibony and the first
author proved in 1995 [\DS] that compact embedded lagrangian submanifolds for a K\"ahler form in
$\C^n$ are rationally convex.

It is natural to try and extend this result to immersed lagrangian submanifolds, which are much
more abundant. The second author obtained in 2000 [\GA]
the rational convexity of generic compact immersed lagrangian submanifolds, those with transversal
self-intersections.

But more complicated self-intersections show up generically when looking at deformations of
lagrangian submanifolds, namely {\it quadratic self-tangencies} : double points where
the two branches are transversal except for one direction along which they get
a quadratic contact. Our result is :

\thm{Theorem} Let $\omega$ a K\"ahler form on $\C^n$, and $L$ a compact immersed lagrangian
submanifold for $\omega$ with quadratic self-tangencies. Then $L$ is rationally convex.
\fthm

As for the former cases, this statement has its counterpart in a symplectic
setting in the spirit of Donaldson's work on symplectic hypersurfaces [\DON] (see [\AGM] for the symplectic version
 of the embedded case). But for simplicity we will stick to our holomorphic context.

\null

The second author would like to thank K. Cieliebak for asking this question. His motivation stemmed from his work on the Floer
homology for pairs of lagrangian submanifolds.

Before entering the proof of our theorem,
 we will review briefly the embedded and the immersed transversal cases.

\null

\subhead 1. The embedded case \endsubhead

\null

Let $L$ be a compact embedded submanifold of dimension $n$ in $\C^n$,
 lagrangian for a K\"ahler form $\omega$.
Choose a potential $\phi$ of this form, i.e. a global function such
 that $dd^c \phi = \omega$. By hypothesis
$d^c\phi$ is a closed form on $L$. By perturbing $\omega$, one can
 always arrange rational periods
for this closed form, or even periods in $2\pi \Z$ after multiplying
 $\omega$ by some constant. Therefore $d^c \phi
=d\psi$  where $\psi$ is a function on $L$ with values in $\R/2\pi \Z$,
and $e^{\phi + i\psi}$ is a well
defined complex-valued function on $L$. It can be extended as a global
 function $h$ on $\C^n$, $\overline {\partial}$-flat
on $L$ (i.e. $\vert \overline {\partial}h\vert =O(d^m)$ for any
 $m$ where $d$
stands for the distance to $L$) and vanishing outside a neighborhood
 of $L$. By construction the first jets of
 $\vert h \vert$ and $e^\phi$ coincide along $L$. Actually it is not hard
 to check that $\vert h \vert_{\phi}\leq e^{-d^2}$ if $d$ is suitably normalized. Here
 $\vert h \vert_{\phi}=\vert h\vert e^{-\phi}$.

Our model for a global
 hypersurface avoiding
$L$ will be $(h^k=0)$. It needs only be corrected by means
 of a $\overline {\partial}$-equation with
 $L^2$-estimates.

Namely consider $f=h^k-v$ where $v$
is a solution of $\overline {\partial} v=\overline {\partial}(h^k)$ satisfying
$\vert \vert v\vert \vert_{2,k\phi} \leq
\vert \vert \overline {\partial}(h^k )\vert \vert_{2,k\phi}$ on a big ball
 (we tend from now on to neglect irrelevant constants).
This translates in $L^{\infty}$-estimates. Indeed we have (cf [\HOR])
$$\vert v(p)\vert \leq 
\epsilon \vert \vert \overline{\partial}v\vert \vert_{\infty, \epsilon}
+ \epsilon^{-n}\vert \vert v \vert \vert_{2,\epsilon}$$
 where the norms are taken on $B(p,\epsilon)$.
Then, forcing the weight, we get  $$\vert v\vert e^{-k\phi}(p)\leq
 e^{k\epsilon} (\epsilon \vert \vert\overline{\partial}(h^k)e^{-k\phi}
\vert \vert_{\infty, \epsilon} + \epsilon^{-n}\vert \vert ve^{-k\phi} 
\vert \vert_{2,\epsilon}).$$ Taking now $\epsilon=\frac{1}{k^2}$ we end up with
$$\vert \vert v\vert \vert_{\infty,k\phi} \leq \frac{1}{k}+k^{2n}\vert \vert v
\vert \vert_{2,k\phi} \leq \frac{1}{k}+k^{2n}\vert \vert \overline {\partial}
(h^k) 
\vert \vert_{2,k\phi}.$$ This last norm is roughly estimated by
$k^{2n+1}(\int_0^1s^{2m}e^{-ks^2}ds)^{\frac {1}{2}}$ using the properties
of $h$. It decays with $k$ if $m$ is large.

Choose now a point $p$ outside $L$. We have $\vert h^k \vert_{k\phi}(p)=0$ but
$\vert h^k \vert_{k\phi}=1$ on $L$. Hence $\vert f \vert_{k\phi}(p) \leq \frac {1}{2}
 \min_L \vert f \vert_{k\phi}$ for $k$ big enough.
Conversely, as in H\" ormander [\HOR] (se also Donaldson [\DON] for peak sections in a symplectic setting), one can construct a global holomorphic function $g$
which picks at $p$ in weighted norm, i.e. satisfying $\vert g \vert_{k\phi}(p) \geq 
 \max_L \vert g \vert_{k\phi}$. For this, take $h$ to be the cut-off of the exponential of the holomorphic 2-jet
of $\phi$ at $p$ and proceed as above.
Therefore $(f(p)g=g(p)f)$ is a global holomorphic
hypersurface passing through $p$ and avoiding $L$. See [\DS] for more details.

\null

\subhead 2. The immersed transversal case \endsubhead

\null

The proof follows exactly the same scheme, the only differences appearing
near the double points. Namely in a neighborhhood of such a point $p$ the function
$h^k$ will be replaced by $h_1^k+h_2^k$, each $h_i$ being constructed as before for each branch $L_i$ of $L$ at $p$. Recall that $\vert h_i\vert_\phi=1$
on $L_i$ and that $\vert h_i\vert_\phi \leq e^{-d_i^2}$ where $d_i$
is the distance to $L_i$.
It remains to guarantee a lower bound on $L$ for $\vert h_1^k+h_2^k \vert_{k\phi}$. Let us treat the case of $L_1$. On this branch we have $\vert h_1^k+h_2^k \vert_{k\phi}
=\vert h_1^k \vert_{k\phi}\vert 1+(\frac {h_2}{h_1})^k\vert=\vert 1+(\frac {h_2}{h_1})^k\vert$. But by construction the first jets of $h_1$ and $h_2$ coincide at $p$.
Hence $\frac {h_2}{h_1}=\frac {h_2e^{-\phi}}{h_1e^{-\phi}}=e^{-d_2^2+iO(d^2)}$ on 
$L_1$ where $d$ is the distance to $p$. Remark now that $d_2$ and $d$ are equivalent on $L_1$ by the transversality of the two branches. Therefore
we get that $\vert h_1^k+h_2^k \vert_{k\phi}=\vert 1+e^{k(-d^2+iO(d^2))}\vert $, and this is bounded from below by the minimum of $\vert 1+e^z \vert$ on a
cone directed by the negative real axis.
Compare with the original proof in [\GA].

\null

\subhead 3. The immersed quadratic case \endsubhead

\null

As before the difficulties appear near the double points. But we cannot rely on such
a simple argument as in the transversal case : the two branches interact too
closely along the quadratic tangency to prevent the vanishing of
 $\vert h_1^k+h_2^k \vert_{k\phi}$ on $L$. We will follow instead the method of
 [\GA].
More precisely, let $L$ be an immersed lagrangian submanifold
for a K\"ahler form $\omega$ with a quadratic
self-tangency at $p$.

\thm{Lemma} There exists a local perturbation of $\omega$
near $p$ such that

i)
$L$ remains lagrangian for the new form

ii) this form has a potential $\phi$ presenting a strict local minimum
0 at $p$

iii) $d^c\phi$ vanishes along $L$ near $p$.

\fthm

Assume this for a moment. The theorem follows with the same method as in the transversal case.
We just replace $h_1^k+h_2^k$ by
$\chi^k+h_1^k+h_2^k$ near $p$ in order to separate the two branches. Here $\chi$ is a cut-off function such that $\chi=e^\epsilon$ ($\epsilon>0$
) near $p$.
It is not hard to check that $\vert \chi^k+h_1^k+h_2^k \vert_{k\phi}$ is bounded from 
below on $L$. Indeed the first term dominates on $L \cap (\phi <\epsilon)$ and one of the other dominates on $L_i \cap (\phi>\epsilon)$.
At the junction $L_i \cap (\phi=\epsilon)$ no cancellation appears :
indeed $h_i$ are real positive along $L_i$ because of the vanishing
of $d^c \phi$.

\thm{Proof of the lemma}\fthm

For simplicity we present it in dimension 2. Take $p$ to be 0
and denote by $z=x+iy, w=u+iv$ the coordinates and by $d$ and
$\delta$ the distances to $0$ and to the $x$-axis respectively.
 These coordinates can be choosen
in such way that $\omega$ is standard at 0, say
$\omega=3dd^c(\vert z\vert ^2+\vert w \vert ^2)+O
(d)$. We may suppose moreover that the tangent planes of the two branches at 0
generate the $(x,u,v)$-hyperplane and intersect along the $x$-axis.

 Each branch is then parametrized by $j:(x,t) \mapsto (x+ip,e^{i\theta}(t+iq)+O(d^3)$, where $\theta$ is a specific angle and $p,q$ are real quadratic
polynomials.

We will construct the potential $\phi$ of the new form
in two steps : firstly we choose a good potential $\psi$
of $\omega$ such that $d^c \psi$ vanishes at a certain order at $0$ along $L$,
namely $d^c \psi=O(d^3+\delta d)$; secondly we slightly perturb $\psi$ in $\phi$ independantly for each branch to achieve the vanishing of $d^c\phi$ along $L$.

We start with $\psi=2x^2+4y^2+3u^2+3v^2+ O(d^3)$. The reason for this particular choice lies in the fact that $d^c(x^2+2y^2)$ vanishes on any
parabola tangent to the $x$-axis at 0 in $\C$. It is not hard to check
that $j^*d^c\psi=\alpha x^2dx+\beta x^2 dt +O(d^3+\delta d)$ where $\alpha$ is the same constant for the two branches, $\beta$ being specific to each branch.
Now we use the freedom to change the rest of order 3 in $\psi$ by pluriharmonic terms. Adding $\gamma \text{Im}(z^3)$ for a suitable real constant $\gamma$ (resp. $\text{Im}(\mu z^2w)$ for some complex number $\mu$)
allows us to kill $\alpha$ (resp. $\beta$ for each branch).

We now perturb $\psi$ to construct $\phi$ near $L_1$ such that $d^c\phi$ vanishes on $L_1$. Without loss of generality we can think of
$L_1$ being flat, say after a rotation the $(x,u)$-plane. 
Indeed it is possible to straighten $L_1$ by means of a local
diffeomorphism tangent to the identity at 0 and $\bar{\partial}$-flat on the branch.
Write $j^*d^c\psi=adx+bdu$ on $L_1$ with $a,b=O(d^3+\delta d)$. Then $j^*d^c(\psi-ay-bv)=0$. 

We then cut this modification off $L_1$ to avoid the interaction with $L_2$. For this, remark that the set $(y^2+v^2 \leq \epsilon (u^2+x^4))$ doesn't intersect $L_2$
except at 0 for a small $\epsilon$. This is a consequence of the quadratic tangency between the two branches. Put $\phi=\psi
-\chi(\frac{y^2+v^2}{u^2+x^4})(ay+bv)$ where $\chi$ is a smooth function with small support on $\R^+$ equal to 1 near 0. It can be checked that $\phi$
is a small $C^2$-perturbation of $\psi$ because of our good control of $a$ and $b$.
Therefore it remains strictly plurisubharmonic and by construction $d^c\phi$ vanishes
on $L_1$. Adding a similar term for $L_2$ concludes.
 \Refs

\widestnumber\no{99}
\refno=0

\bref \by D. Auroux, D. Gayet, J.-P. Mohsen \paper Symplectic hypersurfaces in
the complement of an isotropic submanifold \jour Math. Ann. \vol321\yr2001\pages739--754
\endref
 
\bref \by S. K. Donaldson \paper Symplectic submanifolds and almost-complex geometry \jour J. Diff. Geom. \vol44\yr1996\pages666--705
\endref

\bref \by J. Duval, N. Sibony \paper Polynomial convexity, rational convexity, and currents \jour Duke Math. J. \vol79\yr1995\pages487--513
\endref

\bref \by D. Gayet \paper Convexit\'e rationnelle des sous-vari\'et\'es immerg\'ees
lagrangiennes \jour Ann. Sci. ENS \vol33\yr2000\pages291--300
\endref

\bref \by L. H\"ormander \book The analysis of linear partial differential operators Vol.II \bookinfo Classics in Math.\publ Springer-Verlag \yr 2005 \publaddr Berlin
\endref

\endRefs

\enddocument